\documentclass[12pt,a4]{amsart}

\newtheorem{theorem}{Theorem}[section]

\newtheorem{proposition}[theorem]{Proposition}

\newtheorem*{coro}{Corollary} 
\theoremstyle{definition}

\newcommand{\D}{{\rm I \! D}}
 
\newcommand{\F}{{\mathcal F}}
\newcommand{\g}{{\mathfrak g}}
\newcommand{\HH}{{\mathcal H}}
\newcommand{\CL}{\mathcal {L}}

\newcommand{\CP}{{\mathcal P}}
\newcommand{\pb}{{\mathbf P}}
\newcommand{\R}{{\mathbf R}}

\newcommand{\I}{{\mathcal I}}
\newcommand{\W}{\mathbf W}

\def\paral{/\kern-0.55ex/}
\def\parals_#1{/\kern-0.55ex/_{\!#1}}
\def\ker{\mathop{\rm ker}}

\def\Cyl{\mathop {\rm Cyl}}
\newcommand{\Div}{\operatorname{div}}
\newcommand{\Dom}{\operatorname{Dom}}
\newcommand{\id}{\operatorname{Id}}
\newcommand{\image}{\operatorname{Im}}
\newcommand{\Ric}{\operatorname{Ric^{\#} }}

\newcommand{\pathR}{{\rm I\!R}}

\title[Geometric stochastic analysis on path spaces ]{Geometric stochastic analysis on path spaces }

\author[Elworthy and Li]{K. D. Elworthy and Xue-Mei Li }
\thanks{Xue-Mei Li (Xue-Mei Hairer) has benefited from a Royal Society Leverhulme senior research fellowship. }
\begin{document}

\begin{abstract}
An approach to analysis on path spaces of Riemannian manifolds is described. The spaces are furnished with `Brownian motion' measure which lies on continuous paths, though differentiation is restricted to directions given by tangent paths of finite energy.
An introduction describes the background for paths on $\R^m$ and Malliavin calculus.
For manifold valued paths the approach is to use `It\^o' maps of suitable stochastic differential equations as charts . `Suitability' involves the connection determined by the stochastic differential equation. Some fundamental open problems concerning the calculus and the resulting `Laplacian' are described. A theory for more general diffusion measures is also briefly indicated. The same method is applied as an approach to getting over the fundamental difficulty of defining exterior differentiation as a closed operator, with success for one \& two forms leading to a Hodge -Kodaira operator and decomposition for such forms. Finally there is a brief description of some related results for loop spaces.
\end{abstract}

\maketitle

{\it Classification.}
Primary  58B10; 58J65; Secondary  58A14; 60H07; 60H10; 53C17; 58D20; 58B15.


\section{Introduction}
\subsection{Analysis with Gaussian measures}
Classical differential and geometric analysis is based on Lebesgue measure. The non-existence of an analogue of Lebesgue measure in infinite dimensions is demonstrated by the following theorem:
\begin{theorem}
If $\mu$  is a locally finite Borel measure on a separable Banach space $E$ such that  translations by  every element of $E$ preserve sets of measure zero, then either $\mu=0$ or $E$ is finite dimensional.
\end {theorem}
`Local finiteness' here means that every point of $E$ has a neighbourhood with finite measure. The theorem is a special case of  more general results, e.g see Theorem 17.2 of \cite{Yamasaki-book}. In a sense it is behind many of the mathematical difficulties in 'path integration' and has meant that infinite dimensional  differential and geometric analysis has had to develop its own techniques.

The analysis proposed by Gross was based on his notion of \emph{abstract Wiener spaces}. These are triples $\{i,H,E\}$ where $i:H\rightarrow E$ is a continuous linear
injective map, with dense range, of a Hilbert space $(H, \langle , \rangle_H)$ into a separable Banach space $E$. (Throughout this article \emph{all the linear spaces considered will be real}.) The defining property is that there is a Borel measure $\gamma$, say, on $E$ whose Fourier transform $\hat{\gamma}$ is given by: $$
\hat{\gamma}(l):=\int_E e^{\sqrt{-1}l(x)}d\mu{(x)}=e^{-\frac{|j(l)|_H^2}{2}}$$
for all $l\in E^\ast$ where $j:E^\ast \rightarrow H$ is the adjoint of $i$. This was a generalisation of \emph {classical Wiener space} where some analysis had been previously investigated, particularly by Cameron \& Martin e.g. see \cite {Cameron-Martin}, and also was influenced by Irving Segal's work. It was shown later that all, centred and strictly positive, so-called `Gaussian measures' on a separable Banach space $E$ arise from an essentially unique abstract Wiener space structure on $E$, e.g. see \cite {Dudley-Feldman-LeCamm}.

Classical Wiener space can be considered as the special case when $E$ is the space
$C_0([0,T];\R^m)$ of continuous maps of a fixed interval $[0,T]$ into $\R^m$ which start at the origin, and $H$, sometimes called the \emph {Cameron -Martin space}, is the space of finite energy paths $L^{2,1}_0
([0,T];\R^m)$, i.e. those paths in  $C_0([0,T];\R^m)$ which have distributional derivatives in $L^2$. The map $i$ is the inclusion. The norm for $H$ is given by $ |h|_H ^2=\int_0^t|\dot{h}(t)|^2dt$, and the measure on $E$ is the classical Wiener measure constructed by Wiener, so that the canonical process $[0,T]\times C_0([0,T];\R^m)\rightarrow\R^m$ given by evaluation, is the standard model of Brownian motion. Denote that measure by $ \pb$.

The starting point for Gross's analysis was his extension of Cameron \& Martin's quasi-invariance theorem to the case of abstract Wiener spaces and their measures $\gamma$:

\begin {theorem} Translation by an element $v$  of $E$ preserves sets of measure zero if and only if $v$ lies in the image of $H$. Moreover, for $h\in H$ and integrable $f: E\to \R$ we have for any $t\in \R$: 
\begin{equation}\label{eq:CMform}
\int _E f(x)d\gamma=\int_E f\big(x+ti(h)\big)exp\Big(-t\CP (h)-\frac{t^2}{2}|h|^2_H\Big)d\gamma,
\end{equation}
where $\CP:H\rightarrow L^2(E;\R)$ is the \emph {Paley-Wiener} map.
\end{theorem}
 The Paley-Wiener map is an isometry into $L^2$ defined as the $L^2$-limit of any sequence, $\{l_n\}_{n\ge 1}$ of elements in $E^\ast$ for which $\{j(l_n)\}_{n\ge 1}$ converges in $H$ to $h$. For classical Wiener space it is written as $\sigma \mapsto \int_0^T\langle \dot{h}(s),d\sigma(s)\rangle$ and called the \emph {Paley-Wiener (stochastic) integral}. If $E=H=\R^n$ then $\CP(h)(x)=\langle h,x\rangle_H$ which accounts for the notation $\langle h,x \rangle\tilde{}_H$ sometimes used for it in general.

The Gross-Cameron-Martin formula \label {CMform} was given here with a parameter
$t$ in order to obtain an integration by parts formula from it by differentiating with respect to $t$ at $t=0$. If $f$ is sufficiently regular, for example Fr\'echet differentiable and bounded with bounded derivative $Df:  E\rightarrow E^\ast$,  this yields:
\begin {coro} [Integration by parts]
For all $h\in H$  
\begin{equation}\label{eq:IBP1}
\int_E Df(x)\big(i(h)\big) \;d\gamma(x)
=-\int_E f(x)\Div(h)(x) d\gamma(x) 
\end{equation}
where $\Div(h)=-\CP(h)$.
\end{coro}

 Diffeomorphisms of the form $x\mapsto x+ij\alpha(x)$ between open subsets of $E$
 where $\alpha$ is a $C^1$ map into $E^\ast$ were shown by H-H Kuo to also preserve sets of measure zero. This led him, starting in his thesis, and then others, to the study of \emph{abstract Wiener manifolds}: Banach manifolds modelled on the space $E$ of an abstract Wiener space whose interchange of charts were of Kuo's form, \cite {Kuo-thesis}. Such manifolds have a natural class of Borel measures, locally equivalent to $\gamma$, and many of the usual constructions and results of the finite dimensional situation go over to them, \cite {Kuo-diffusions}, \cite {Elworthy-Trieste}, \cite{Eells-Elworthy-CIME}, \cite {Piech-82}. 
 
 For any abstract Wiener space the map $i$ is compact. It follows that the derivatives
 of transformations of Kuo's type are linear maps which differ from the identity by a compact operator. This implies that abstract Wiener manifolds are Fredholm manifolds, \cite{Eells-Elworthy-CIME}. For a wide class of Banach spaces $E$ the theory and classification of such manifold structures showed that every separable metrisable manifold $M$ modelled on $E$ can be given the structure of an abstract  Wiener manifold, with the K-theory of $M$ playing the major role in their classification, \cite{Eells-Elworthy-ICM}. 
 
 It soon became clear that although interesting manifolds, such as path and loop spaces on finite dimensional manifolds admit these structures, in most interesting cases there is no natural one. Exceptions are finite codimensional submanifolds of abstract Wiener spaces, such as the space of paths from one submanifold embedded in $\R^m$
 to another. Also see \cite {Freed}. More general transformations preserving sets of measure zero were described, notably by Ramer in 1974, and then using Malliavin calculus by Kusuoka in 1982,
 and for flows of a class of vector fields on classical Wiener space by Driver, \cite{Driver92}. However the form of these transformations is not so different from those of Kuo, though the identity map in the decomposition may be replaced by a  `rotation'. See \cite {Ustunel-Zakai-book}. This together with the advent of Malliavin calculus, in 1976, with emphasis on mappings determined by stochastic differential equations, led to a move away from this approach, or at least a major modification of it \cite {Kusuoka-88}, \cite {Kusuoka-ICM}. 
 
 In  \cite {Gross-AWS} Gross shows that for any abstract Wiener space $\{i,H,E\}$ there is an abstract Wiener space $\{i',H,E'\}$ and a compact linear map $k:E'\rightarrow E$ such that  $i=k\circ i'$. In other words, in the infinite dimensional case the
 measure $\gamma$ can be considered to lie in a smaller space than $E$; (however $H$ itself has measure zero). In the classical case this is demonstrated by the fact that the space of continuous functions can be replaced by the closure of $L^{2,1}_0$ in the space of H\"older continuous functions of exponent $\alpha$ for any $0<\alpha<1/2$. In Malliavin calculus on these linear spaces the space $E$ loses its importance, and in some treatments essentially disappears, e.g.  see \cite{Ito-Charingworth}, and \cite {Malliavin-book}.
In the latter it is the Paley-Wiener functions, $\{\CP (h):h\in H\}$ which play the dominant role, returning to Segal's `weak distribution' theory, \cite {Segal-65}. However
in the non-linear case of diffusion measures on path spaces of manifolds  it seems necessary, at least at the moment, to deal with the actual manifold on which the measures sit, though this could be taken to be H\"older continuous paths rather than continuous paths if that is more convenient. The treatment of Malliavin calculus below is organised with this in mind.

\subsection {Malliavin calculus on $E$}
From Gross's work, especially \cite{Gross-potential}, it was clear that the basic differentiation operator on an abstract Wiener space should be the \emph{$H$-derivative}. This could be defined on a basic domain, $\Dom (d_H)$, of functions $f:E\rightarrow\R$ consisting of a set of Fr\'echet differentiable functions which is dense in $L^p$ and whose $H$-derivatives: $ d_Hf: E\rightarrow H^*$ given by $d_Hf_x(h)=Df(x)(i(h))$ lie in  $L^p$,  for all $1\le  p<\infty$. The integration by parts formula,  equation (\ref{eq:IBP1}), implies that $d_H$
is closable as a map between $L^p$ spaces with closure a closed linear map $$ d:\Dom(d)\subset L^p(E;\R)\rightarrow L^p(E;H^\ast).$$ 
Let $\D^{p,1}$ denote $\Dom(d)$ with its graph norm.

Our Paley-Wiener functionals, $\CP (h)$,  are easily seen to be in $\D^{p,1}$
for all $1\le  p<\infty$ with $d \CP (h)_x(k)=\langle h,k\rangle_H$ for  all $x\in E$ and $k\in H$, despite their lack of continuity in $E$. In fact the main point of the theory is that, for classical Wiener space, more general stochastic integrals and solution maps of stochastic differential equations, as described below, all lie in these Sobolev
spaces.

The following characterisation of $\D ^{p,1}$ for $1<p<\infty$ was given by Sugita:
\begin {theorem}[\cite{Sugita}]
 If $f\in L^p(E;\R)$ then $f\in\D ^{p,1}$ if and only if both of the following hold
\begin{enumerate}
\item For each $h\in H$ there is a function $f_h:E \times\R\rightarrow\R$ which is absolutely continuous in the second variable and has $f_h(x,t)=f(x+th)$, for almost all $x\in E$, for each $t\in\R$.
\item There exists $df\in L^p(E;H^\ast)$ such that for any $h\in H$, $\frac{1}{t}(f(x+th)-f(x))$ converges in measure to $df_x(h)$ as $t\rightarrow 0$.
\end{enumerate}
\end{theorem}
From this we see that the spaces $\D ^{p,1}$, for $1<p<\infty$ are independent of any reasonable choice of initial domain $\Dom(d_H)$. A comforting fact; but one which is still open in the corresponding situation for paths on curved spaces, as will be described below.

For functions with values in a separable Hilbert space $G$ the spaces $\D ^{p,1}(E;G)$ are defined in the analogous way,  with the derivative 
$df$ now mapping $E$ into $\mathcal{L}_2(H;G)$, the space of Hilbert-Schmidt operators of $H$ into $G$. This is a Hilbert space often identified with the completed tensor product $G\otimes H$. It occurs because a basic property of an abstract Wiener space is that any continuous linear map from $E$ to a Hilbert space $G$,
such as $Df(x)$ if $f:E \rightarrow G$ is Fr\'echet differentiable, gives a Hilbert-Schmidt operator when composed with $i$, e.g. see Thm 17.3 in \cite{Yamasaki-book}. 
 We can iterate this procedure to obtain higher order Sobolev spaces.

As usual the gradient can be defined for functions in $\D ^{p,1},1<p<\infty$, by $\langle\nabla f(x), h \rangle_H=df_x(h)$ to give an \emph{$H$-vector field}, $\nabla f:E\rightarrow H$. It is a closed operator between the $L^p$ spaces with the negative of its adjoint denoted by $\Div$, a closed operator from $\Dom(\Div)$ in $L^q(E; H)$ to $L^q(E;\R)$, where $q$ is the conjugate of $p$. Similar we have the adjoint $d^*$ of $d$.
From this we get the analogue of the (Witten) Laplacian, or the `Ornstein-Uhlenbeck operator', $\mathcal{L}= \Div \nabla=-d^*d$. In the case $E=H=\R^n$
this is given by:$$ \mathcal{L}(f)(x)= \triangle f(x)-Df(x)(x)$$ for $\triangle$ the usual Laplacian (with negative spectrum) of $\mathbb{R}^n$.

The Ornstein -Uhlenbeck operator acting on $L^2$ is the well known operator whose spectrum consists of $0$ as unique ground state, together with the negative integers as eigenvalues of infinite multiplicity, corresponding to the homogeneos chaos decomposition of $L^2(E;\R)$, and conjugate to the number operator of mathematical physics acting on the real symmetric Fock space. For example from above we see that for $h\in H$ the map $\CP (h)$ is an eigenvector of eigenvalue minus one (so giving the `one-particle' space). For more, see for example \cite{Hsu-book}, \cite{Nualart-book}, or \cite{Ikeda-Watanabe}.

For classical Wiener space an H-vector field $V:C_0\rightarrow L^{2,1}_0$ is said to be
\emph{non-anticipating} if for each time t its value $V(\sigma)_t$ at the path $\sigma$
depends only on the restriction of $\sigma$ to the interval $[0,T]$. If this holds and it is in 
$L^2$, then it is in the domain of the divergence operator and $\Div(V)(\sigma)$
is  precisely the negative of the \emph{It\^o stochastic integral}, $\int_0^T\dot{V}(\sigma)_t\; d\sigma(t)$,
as shown by Gaveau. This is the integral which is the basis of stochastic calculus. In the  
anticipating case it is the \emph {Skorohod}, or \emph{Ramer-Skorohod integral}, now by  definition: although here the word `integral' can be misleading since, as in finite dimensions, differentiation may be involved, \cite {Nualart-book}.

An $L^2$-deRham and Hodge-Kodaira theory was given in this context by Shigekawa \cite{Shigekawa86}. The k- forms were `$H$-forms',  i.e. maps from $E$ into  $\wedge^kH$ where $\wedge^k$ denotes the Hilbert space completion of the k-th exterior power, with the exterior derivative being a closed operator derived from our $H$-derivative $d$. The Hodge decomposition was just as in finite dimensional, standard, $L^2$-theory, and Shigekawa proved a vanishing theorem, implying the expected triviality of the deRham cohomology.  A theory of finite co-dimensional forms was proposed by Ramer in his Thesis, in the context of abstract Wiener manifolds; further developments were made by Kusuoka, \cite{Kusuoka-foundationsII}, but more is needed  to develop the theory, even on domains in these linear spaces.

\section{Scalar analysis on paths in $M$}
\subsection{Brownian motion measure and Bismut tangent spaces}\label{sec-Bismut}
Consider a smooth manifold $M$. For a fixed time $T>0$, and a fixed $x_0\in M$ let 
$C_{x_0}([0,T];M)$, or simply $C_{x_0}$,  denote the space of continuous paths $\sigma :[0,T]\rightarrow M$
starting at $x_0$, together with its usual $C^\infty$  Banach manifold structure, e.g.see \cite{Eliasson67} or \cite{Palais}. The tangent space $T_\sigma C_{x_0}$ to $C_{x_0}$ at a point $\sigma$ can be identified with the space of continuous paths $v:[0,T]\rightarrow TM$ into the tangent bundle to $M$, such that $v(0)=0$ and $v(t)\in T_{\sigma(t)}M$ for $0\le t\le T$.
For a complete Riemannian manifold the \emph{ Brownian motion measure}, $\mu_{x_0}$,
on $C_{x_0}$ is the unique Borel measure for which 
\begin{equation}\label{eq:BM}
\mu_{x_0}\Big(\{\sigma\in C_{x_0}:\sigma(t_j)\in A_j, j=1,2,\dots,k\}\Big)
 =\int_{A_1}\int_{A_2}\dots \int_{A_n}\prod_{j=0}^{j=k-1} p_{t_{j+1}-t_j}(x_j;dx_{j+1})
\end{equation}
where $0=t_0<t_1<\dots<t_k\le  T$, the $A_j$ are Borel subsets of $M$, and
the measures $p_t(x,dy)$ are the heat kernel measures: $p_t(x,dy)=p_t(x,y)dy$ for $p_t(x,y)$ the fundamental solution of the heat equation $\frac{\partial f}{\partial t}=\frac{1}{2}\triangle$ for $\triangle$ the Laplace Beltrami operator, div grad, of $M$.

For simplicity we shall assume that $M$ is compact. Let its dimension be $n$.

From the successes of the flat space case it was expected that, to do analysis on the path space $C_{x_0}$ using Brownian motion measure, the differentiation should only take place in a special set of directions. In the case of Gaussian measures on linear space a natural choice was given, as described above, by the linear structure together with the measure: but there are other choices as we see in section \ref{sec-general} below  and it is not clear if the measure plus the differential structure does determine a special one, c.f. \cite {Elworthy-Dublin}. Nevertheless a natural choice for Brownian motion  measure is the \emph {Bismut tangent spaces}. These are Hilbert spaces, $\mathcal{H}_\sigma, $ of tangent vectors, defined for almost all $\sigma \in C_{x_0}$ by 
\begin {equation} \label{eq:Bismut}
\mathcal{H}_\sigma=\{v\in T_\sigma C_{x_0}: (\parals_\cdot)^{-1}v(\cdot)\in L^{2,1}([0,T];T_{x_0}M)\}
\end{equation}
where $\parals_t$ denotes parallel translation along $\sigma$ using the Levi-Civita connection.

Because our paths $\sigma$ are typically so irregular, e.g. almost surely $\alpha$-Holder continuous only for $\alpha <1/2$, the parallel translation has to be constructed by stochastic differential equations and so is only defined along almost all paths. However if we set $\mathcal{H}=\cup_\sigma\mathcal{H}_\sigma\subset TC_{x_0}$ we will see that it has the rudiments of a bundle structure.  We will call its sections $H$-vector fields, and the sections of its dual bundle $\mathcal{H}^*$ will be called $H$-one-forms, c.f. \cite {Jones-Leandre}.

In \cite{Driver92}, Driver extended Cameron-Martin's theorem and the formulae (\ref{eq:CMform}), and (\ref{eq:IBP1}) to this situation, showing that if $V^h$ is the $H$-vector field
whose value at $\sigma$ is obtained by parallel translation of a fixed element $h\in L^{2,1}([0,T];T_{x_0}M)$ along $\sigma$ then this measurable vector field has a solution flow which preserves sets of $\mu_{x_0}$-zero, with consequent analogues of equations
(\ref{eq:CMform}) and (\ref{eq:IBP1}).

As for flat space the integration by parts formula gives closability of the $H$-derivative
$d_H:\Dom(d_H)\to L^2\Gamma\mathcal{H}^\ast$ from its domain in $L^2$ into the $L^2$ -$H$-one-forms. It works for the $ L^p$-spaces but we shall only mention $L^2$ from now on for simplicity.  A natural, essentially the smallest natural, domain to choose is to let $\Dom(d_H)$ be the space $\Cyl(M)$ of smooth cylinder functions: those maps of the form $\sigma \mapsto F(\sigma(t_i),\dots\sigma(t_k))$ for some smooth $F$ defined on the $k$-fold product of $M$, some $0\le t_1<\dots<t_k\le  T$, any natural number $k$. Other choices include the space of (Fr\'echet) $C^1$- functions which are bounded together with their derivatives, using the natural Finsler metric on $C_{x_0}$. However this time we do not know if these lead to the same domain for the closure of $d_H$, see \cite{Elworthy-Li-Trento, Elworthy-Li-Acta-preprint}.

We must make a choice, and will choose $\Cyl(M)$ as basic domain. With this choice let
$\D ^{2,1}(C_{x_0})$, or $\D ^{2,1}$, denote the domain of the $L^2$-closure of $d_H$ with its graph norm, with  $\D ^{2,1}(C_{x_0};G)$ for the corresponding space of $G$-valued functions, $G$ a separable Hilbert space. Let $d$ denote the closure of $d_H$, so if $f: C_{x_0}\rightarrow G $ is in $\D ^{2,1}(C_{x_0};G)$ then $df$ is an $L^2$-section of $\mathcal{L}_2(\mathcal{H};G)$ the `bundle' with fibre at $\sigma$ the space of Hilbert-Schmidt maps of $\mathcal{H}_\sigma$ into $G$, sometimes denoted by $\mathcal{G}\otimes\mathcal{H}$.

With this we get a closed operator $\nabla$ as usual, mapping its domain $\D ^{2,1}$ into $H$-vector fields, with adjoint the negative of a closed operator $\Div$. As usual we have a self adjoint `Laplacian', or `Ornstein-Uhlenbeck' operator
$\mathcal{L}$ acting on functions, defined by $\mathcal{L}=\Div\nabla=-d^*d$. The associated Dirichlet forms and processes have been studied, e.g. \cite {Driver-Rockner}, \cite{Eberle-book}. Norris devised a stochastic partial differential equation to construct associated `Brownian motions' or `Ornstein-Uhlenbeck processes' on these path spaces, treating them as two parameter $M$-valued processes, \cite{Norris-pde}. The existence of a spectral gap for $\CL$ was proved by S. Fang,  and Log Sobolev inequalities independently by E. Hsu and Aida \& Elworthy, see \cite {Hsu-book}, \cite{Elworthy-LeJan-Li-book}. However little, if anything, appears to be known else about its spectrum.

To discuss higher derivatives it is convenient to have a `connection' 
on $\mathcal{H}$ in order to differentiate its sections. The most obvious choice is to use the trivialisation of $\mathcal{H}$ obtained simply by parallel translating every element in 
each $\mathcal{H}_\sigma$ back to an element of  
$L^{2,1}([0,T];T_{x_0}M)$, so that $H$-vector fields can be considered as maps of $C_{x_0}$ into $L^{2,1}([0,T];T_{x_0}M)$
to which we may try to apply our closed derivative operator $d$. 
This approach was used effectively, for example in \cite{Leandre-96}. However it does not conserve the $C_{\id}([0,T];GL(n))$-structure 
of our path space, nor as Cruziero\& Malliavin pointed out, does it 
fit well with the underlying `Markovianity' of our set up. This led 
them to
the `Markovian' connection, see \cite{Cruzeiro-Malliavin-96},
a modification of which we will describe below. 

\subsection{It\^o maps and the stochastic development} 
\label {sec-Ito.maps}

The stochastic development map $\mathcal{D}:C_0([0,T];T_{x_0}M)\rightarrow C_{x_0}$ is an almost surely defined version of the Cartan development, describing `rolling without slipping' along smooth paths. Its inverse is given by $\mathcal{D}^{-1}(\sigma)(t)=\int_0^t (\parals_t)^{-1}\circ d\sigma(t)$, where the integral is a Stratonovitch stochastic integral, and $\parals_t$ refers to parallel translation along the path $\sigma$, (defined for almost all paths).  Reformulating Gangolli, \cite {Gangolli},\cite {Eells-Elworthy-CIME}, it was shown by Eells \& Elworthy that it sends Wiener measure to the Brownian motion measure. A fundamental result of Malliavin calculus is that, for each time $t$ the map can be H-differentiated infinitely often in the Sobolev sense. This was used by Driver to transfer
his results about flows of vector fields, and integration by parts formulae, from flat space
to $C_{x_0}$, see \cite{Driver92} where background details are included. However
the use of $\mathcal{D}$ as a chart was limited because its $H$-derivative does not map $L^{2,1}([0,T];T_{x_0}M)$ to the Bismut tangent spaces. Furthermore from  \cite {X-D-Li-2003}
it now seems that, unless $M$ is flat, composition with $\mathcal{D}$ will not pull elements in $\D^{2,1}(C_{x_0})$ back to elements in the domain of $d$: there will be a loss of differentiability.

An alternative technique is to use the solution maps, \emph{It\^o maps}, of more simple stochastic differential equations as replacements for charts. For this take a (Stratonovich) stochastic differential equation \begin {equation} \label{eq:SDE}
 dx_t=X(x_t)\circ dB_t +A(x_t) dt \end{equation}
 on $M$. Here $A$ is a smooth vector field and $X$ gives linear maps $X(x):\R^m \rightarrow T_xM$,  smooth in $x\in M$. Also $B$ is the canonical Brownian motion given by $B_t: C_0([0,T];\R^m)\rightarrow \R^m$ with $B_t (\omega)=\omega(t)$ for $C_0([0,T];\R^m)$ furnished with its Wiener measure, which we shall now denote by $\pb$.
 
The solution $x_t:C_0([0,T];\R^m)\rightarrow M$ to such an equation, starting from $x_0$, can be obtained by `Wong-Zakai approximation': taking piecewise linear approximations $B^\Pi_t$ to the Brownian motion, for each partiition $\Pi$ of $[0,T]$, and solving the family of ordinary differential equations $$\frac{dx^\Pi(\omega)}{dt}=X(x^\Pi_t(\omega))\frac{dB^\Pi (\omega)}{dt} +A(x^\Pi_t(\omega))$$ starting at $x_0$, for each $\omega$. The required solution $x_t$ is given by $x_t(\omega)=\mathcal{I}(\omega)_t$ for $\mathcal{I}$ the limit in probability of $ x^\Pi: C_0([0,T];\R^m)\rightarrow C_{x_0}$ as the mesh of $\Pi$ goes to zero. The map  $\mathcal{I}$ is the \emph{ It\^o map}. To be precise we have to choose it as a  representative from an almost sure equivalence class of measurable maps. However, as with the stochastic development these maps can be differentiated arbitrarily many times in the sense of Malliavin calculus. In particular for almost all $\omega$ there is a linear $H$-derivative $T_\omega\mathcal{I}: H\rightarrow T_{\mathcal{I}(\omega)}C_{x_0}$.

The solutions to equation~(\ref{eq:SDE}) form a Markov process with generator $\mathcal{A}$ where \begin {equation}  \label{eq:Horform}
 \mathcal{A}={1\over 2} \sum_{j=1}^{m}\mathcal{L}_{X^j}\mathcal{L}_{X^j} +\mathcal{L}_A. \end{equation}
For them to be Brownian motions we need $\mathcal{A}={1\over 2}\triangle$ which requires  each $X(x):\R^m\rightarrow T_xM$ to be surjective and induce the given Riemannian metric on the tangent space, or equivalently for the adjoint $Y_x: T_xM\rightarrow \R^m$ of $X(x)$ to be a right inverse of $X(x)$. Given that, we may choose the vector field $A$ appropriately. Then $\mathcal{I}$ will map the flat Wiener measure $\pb$ to our Brownian motion measure $\mu_{x_0}$.
 In general the dimension, $m$, of the space on which the driving Brownian motion runs, will be larger than that of $M$ so that $\mathcal{I}$ will not be injective. The disadvantage of this can be reduced by `filtering out the redundant noise' and to do this successfully we need to note that our SDE  for Brownian motion determines a metric connection,$\breve{\nabla}$ on $TM$ by using $X$ to project the trivial connection on the trivial $\R^m$-bundle onto $TM$: for a vector field $U$ and tangent vector $v\in T_xM$ the covariant derivative of $U$ in the direction $v$ is given by 
 \begin{equation}\label{eq:LJW} 
 \breve{\nabla}_vU=X(x)\Big(d[y\mapsto Y_yU(y)]_x(v)\Big).
  \end{equation}
   It follows from Narasimhan \& Ramanan's theory of universal connections
 that every metric connection on $TM$ can be obtained by a suitable choice of $X$, see \cite{Elworthy-LeJan-Li-book}, or \cite {Quillen} for a direct proof. To obtain the Levi-Civita connection we can use Nash's theorem to take an isometric embedding $j:M\rightarrow\R^m$ for some $m$ and then set $X(x)=(dj)^\ast_x$, the adjoint of
 $(dj)_x$. With $A=0$ the resulting `gradient' SDE has Brownian motions as solutions as required. For Riemannian symmetric spaces it may be useful to use the homogeneous space structure; for example if $M$ is a compact Lie group with bi-invariant metric we may take $\R^m$ to be a copy of the direct sum $\g\oplus\g$  of the Lie algebra, $\g=T_{\id}M$, of $M$ with itself and define $X(x)(e,e')=TR_x(e)-TL_x(e')$ with $A=0$,  where $TR_x$ and $TL_x$ are the derivatives of left and right translation by $x$,  \cite{Elworthy-LeJan-Li-book}.
 
   One basic result, extending estimates in \cite{Aida-Elworthy}, which contrasts with the stochastic development is:
 \begin{theorem}
 [\cite{Elworthy-Li-Acta-preprint}]
 Suppose the connection $\breve{\nabla}$ induced by the SDE is the Levi-Civita connection. The the pull back by $\I$ of cylindrical one-forms on $C_{x_0}$ extends to a continuous linear map $\I^*: L^2\mathcal{H}^\ast \to L^2(C_0([0,T];\R^m);H^*)$ of $L^2$ $H$-one-forms on $C_{x_0}$ to those on the flat path space.
 \end{theorem}
 
Here for a cylindrical, or other one form, $\phi$,  on $C_{x_0}$,  the pull-back $H$-form $\mathcal{I}^*(\phi)$ is given by $\mathcal{I}^\ast(\phi)_\omega(h)=\phi(T_\omega\mathcal{I}(h))$ for $h\in H$. However, in general the $H$-derivative $T\mathcal{I}$ does not map $H$ into the Bismut tangent spaces and so for $H$-one-forms $\phi$ the pullback does not have a classical meaning, though it does have an expression as an It\^o integral under our condition on $\breve{\nabla}$. If $\breve{\nabla}$ were not the Levi-Civita connection this integral would be a Skorohod integral with a consequent loss of differentiability expected, as for the stochastic development map in \cite {X-D-Li-2003}.
There is an important equivalent dual, or `co-joint', version to this result. For this suppose
$\alpha:C_{x_0}\rightarrow H$ is an $H$-vector field in $L^2$. For almost all $\sigma\in C_{x_0} $ we can `integrate over the fibre of $\mathcal{I}$' at $\sigma$ to obtain $\overline{T\mathcal{I}(\alpha)}_\sigma\in T_\sigma C_{x_0}$. Mathematically that is achieved by taking the conditional expectation with respect to the $\sigma$-algebra $\F^{x_0}$ on $C_{x_0}$ generated by $\mathcal{I}$: 
$$\overline{T\mathcal{I}(\alpha)}_\sigma=\mathbb{E}\big\{T\mathcal{I}(\alpha(-)) \; |\; \mathcal{I}(-)=\sigma\big\}.$$
\begin {theorem} 
[\cite{Elworthy-Li-Acta-preprint}]
\label{th:TI}   Suppose the connection $\breve{\nabla}$ induced by the SDE is the Levi-Civita connection. Then for all $H$-vector fields $\alpha$ in  $L^2$, we have $\overline{T\mathcal{I}(\alpha)}_\sigma\in \mathcal{H}_\sigma$ almost surely, giving a continuous linear map $\overline{T\mathcal{I}(-)}: L^2(C_{x_0};H)\rightarrow L^2\mathcal{H}$.
 \end{theorem}
 
 When $\alpha$ is constant, with value $h$ say, we write $\overline{TI}_\sigma(h)$ for $\overline{T\mathcal{I}(\alpha)}_\sigma$. This map was known earlier, \cite{Elworthy-LeJan-Li-book}, to map $H$ isomorphically onto $ \mathcal{H}_\sigma$, with our assumption on $\breve{\nabla}$. In fact it has the explicit expression:
 \begin {equation} \label{eq:TI}
  \overline{TI}_\sigma(h)_t=\W(X(\sigma(-))\dot{h})
 \end {equation} where $\W : L^2([0,T];TM)\to \mathcal{H}$ is an isomorphism of the Bismut tangent `bundle', where defined, with the $L^2$-tangent bundle $L^2TC_{x_0}$ of $C_{x_0}$:
 \begin{eqnarray*}
 L^2T_\sigma C_{x_0}
 =&\{v:[0,T]\rightarrow TM  \hbox{ such that}  &v(t)\in T_{\sigma(t)}M,\;
  0\le t\le T \\
 &&\hbox{ and }  \int_0^T|v(t)|_{\sigma(t)}^2dt <\infty \}. \\
\end{eqnarray*}
 The isomorphism is the inverse of the `damped derivative' along the paths of $C_{x_0}$ 
 \begin{equation}
 \label {eq:damp}
 \frac{\D}{dt}=\frac{D}{dt}+{1\over 2} \Ric: \mathcal{H}\rightarrow L^2TC_{x_0}
 \end{equation}
 where $\Ric:TM\rightarrow TM$ corresponds to the Ricci curvature.
 
 It is convenient to give $\mathcal{H}$ the Riemannian metric and bundle structure it inherits from this isomorphism with the bundle of $L^2$ `tangent vectors'. The latter is a smooth Hilbert bundle over $C_{x_0}$ with structure group $C_{\id}([0,T];O(n))$. It also has a natural metric, `Levi-Civita',  connection, the `pointwise connection' induced from the Levi-Civita connection on $M$, \cite {Eliasson67}. Moving this to $\mathcal{H}$ by $\W $ gives a metric connection which is easily seen to be that projected onto $\mathcal{H}$ by $\overline{T\mathcal{I}}$, in the same way as we defined $\breve{\nabla}$. This connection agrees with the `damped Markovian' connection of 
 Cruzeiro \& Fang, see \cite{Cruzeiro-Fang-97}, referred to above. It can be used to define higher order derivative  operators and Sobolev spaces, and Sobolev spaces of sections of $\mathcal{H}$, e.g. $\D ^{2,1}\mathcal{H}$,  the domain of the $L^2$-closure
 of the covariant $H$-derivative acting on sections of $\HH$. The latter is shown to be in the domain of $\Div$ in \cite {Elworthy-Li-Acta-preprint}: a result proved by M. \& P. Kree for classical Wiener measure in 1983.
 
We can define an $L^2$-function $f:C_{x_0}\rightarrow\R$ to be \emph{weakly differentiable} if it is in the domain of the adjoint of the restriction of $\Div$ to  $\D ^{2,1}\mathcal{H}$. Let $W^{2,1}$ denote the space of such functions with its graph norm. Thus for $f\in W^{2,1}$ there exists $\widetilde{df}\in L^2\mathcal{H^\ast}$ such that if $V\in \D^{2,1}\mathcal{H}$ then $$
\int_{C_{x_0}}f(\sigma)\Div(V(\sigma)) d\mu_{x_0}=-\int_{C_{x_0}}\widetilde{df}(V(\sigma)) d\mu_{x_0}.$$
For paths on $\R^m$ it follows from \cite {Sugita} that weak differentiability implies differentiability, in our Sobolev sense. 

We have the following intertwining result:
\begin{theorem}[\cite {Elworthy-Li-Acta-preprint}]
\label{th:intertwining}
Suppose the connection $\breve{\nabla}$ induced by the SDE is the Levi-Civita connection. Then $f\in W^{2,1}$ if and only if $f\circ\mathcal{I}\in \D ^{2,1}C_{0}([0,T]; \R^m)$ and composition with $\mathcal{I}$ gives a continuous linear map of $W^{2,1}$ onto the space $ \D ^{2,1}_{\mathcal{F}^{x_0}}$ of those elements in $\D ^{2,1}C_{0}([0,T];\R^m)$ which are $\mathcal{F}^{x_0}$-measurable. Moreover for $f\in W^{2,1}$ we have 
$$d(f\circ\mathcal{I})=\mathcal{I}^\ast \widetilde{df}.$$
\end{theorem}
Preliminary versions of some of the above results were given in \cite {Elworthy-Li-Hodge-1}. 
A fundamental question is whether $W^{2,1}=\D ^{2,1}$. Applying results of Eberle, \cite{Eberle-book}, it is shown in \cite {Elworthy-Li-Acta-preprint} that this equality holds if and only if \emph{Markov uniqueness} holds for the operator $\mathcal{L}$ defined above but with domain $\Cyl(M)$. Markov uniqueness is a weaker notion than essential self-adjointness. Probabilistically it relates to uniqueness of solutions to the martingale problem, and it essentially means that there is a unique extension which generates a Markov semigroup. Equality would also imply the independence of $\D ^{2,1}$ from the choice of initial domain $\Dom(d_H)$. We do not know of any non-flat manifolds $M$ for which an answer is known to these questions.
A positive answer would follow from a positive answer to the following:
\begin {itemize} 
\item
If $f\in  \D ^{2,1}C_{0}([0,T];\R^m)$ is its conditional expectation $\mathbb{E}\{f|\mathcal{F}^{x_0}\}$ also in $\D ^{2,1}?$
\end{itemize}
This is described concisely in \cite{Elworthy-Li-Trento}, and in detail in \cite{Elworthy-Li-Acta-preprint}, describing some partial results and correcting claims made in our 2004 Comptes-Rendues note. A discussion somewhat related to the above question, by Airault, Malliavin \& Ren, is in \cite{Airault-Malliavin-Ren-04}.

\subsection{More general diffusion measures}
\label{sec-general}
Let $\mathcal{A}$ be a smooth diffusion generator on $M$ i.e. it is a semi-elliptic second order differential operator with no zero order term, acting on real valued functions on $M$. Essentially as for the case $\mathcal{A}=\frac{1}{2}\triangle$, there is an induced measure $\mu^\mathcal{A}_{x_0}$ on $C_{x_0}$.
 
 To extend the previous results to  do analysis with such a measure we will suppose the principal symbol $\sigma^\mathcal{A}: T^\ast M\rightarrow TM$  of $\mathcal{A}$ has constant rank, and so has image in a sub-bundle $E$ of $TM$. This is equivalent to requiring that $\mathcal{A}$ has a Hormander form, as equation (\ref{eq:Horform}), with the vector fields $X^j$ being sections of $E$.   
 
 In general there is now no obvious choice of a connection with which to define `Bismut 
 tangent' spaces. We therefore choose any metric connection on $E$ and as before, using Narasimhan \& Ramanan's theorem, take a stochastic differential equation (\ref{eq:SDE}) for which the induced connection $\breve{\nabla}$ on $E$ is that chosen one, and for which (\ref{eq:Horform}) holds. To define the the Bismut tangent spaces it is convenient to use the \emph{adjoint semi-connection}, $\widehat{\nabla}$, which allows differentiation of all smooth vector fields, but only in $E$-directions. It is defined by 
 $$\widehat{\nabla}_{U(x)}V=\breve{\nabla}_{V(x)}U+[U,V](x)\in T_xM$$
 for $U$ a section of $E$ and $V$ a vector field on $M$, \cite {Elworthy-LeJan-Li-book}.
 
 Adjoint connections were used in a similar way in order to use different Bismut tangent spaces for Brownian measures, by Driver in \cite{Driver92}. The adjoint of the Levi-civita connection is itself; that of the flat left invariant connection  on a Lie group is a  flat right invariant connection. For more examples see \cite{Elworthy-LeJan-Li-book}.  Semi-connections are also called `partial connections' or `$E$-connections'.
 
 We now define $\mathcal{H}_\sigma$ to be the set of those $v\in T_\sigma C_{x_0}$ for which $\frac{\hat\D} {dt}(v)\in L^2([0,T];E)$ where \begin{equation}  
  \frac{\hat\D} {dt}=\frac{\hat D}{dt}+\frac{1}{2} \breve{\Ric}-\breve{\nabla}_-A 
   \end{equation}
 where the covariant differentiation is done using the semi-connection while $\breve{\Ric}:TM\rightarrow E $ corresponds to the Ricci curvature for $\breve{\nabla}$.
 If $A$ does not take values in $E$ then this operator needs special interpretation, \cite{Elworthy-LeJan-Li-book}. Since $L^2([0,T];E)\bigcap L^2TC_{x_0}$ is a smooth Hilbert bundle, as for the case $E=TM$, with pointwise connection induced from $\breve{\nabla}$, we can induce all this structure, at least almost surely, on $\mathcal{H}$. When there is a metric on $TM$ to which  the semi-connection, $\widehat{\nabla}$, is adapted the theory goes essentially as before, \cite {Elworthy-Li-Acta-preprint}. If not there may be some loss of integrability in the intertwining, for example, but the operator $\mathcal{L}$
has a spectral gap; indeed there is a Log Sobolev inequality, \cite{Elworthy-LeJan-Li-book}. The Dirichlet forms which arise in this situation are discussed in \cite {Elworthy-Ma}.
 
\section{Towards an $L^2$ deRham-Hodge-Kodaira Theory}
\subsection{The spaces of $H$-forms}
Following Shigekawa's rather complete  $L^2$- deRham theory for $H$-forms on abstract Wiener spaces it would be natural to base such a theory on sections of the dual bundles to the exterior products $\wedge^k\mathcal{H}$ of the Bismut tangent bundle, using the Hilbert space completion of the exterior powers of each $\mathcal{H}_\sigma$. 
However this runs into difficulties even at defining the exterior derivative of an $H$ one-form, $\phi$, say: Recall that the standard formula for the exterior derivative $d\phi$ is 
$$d\phi\big(U(x)\wedge V(x)\big)
=\mathcal{L}_U\big(\phi(V(-))\big)(x)-\mathcal{L}_V\big(\phi(U(-))\big)(x)-\phi\big([U,V](x)\big)$$ for vector fields $U$ and $V$.
However if $U$ and $V$ are $H$-vector fields their bracket need not be and so if      $\phi$ is an $H$-form the last term in the expression above will not in general be defined.
One way round this is to interpret this final term as a stochastic integral, in general a Skorohod integral. This was carried through by 
L\'eandre in \cite{Leandre-96} where he obtained a deRham complex in this situation and for loop spaces, proving that the resulting deRham  cohomology agrees with the topological real cohomology. However
this was not really an $L^2$ theory and did not include a version of the Hodge-Kodaira  
Laplacian.

A proposal made in \cite{Elworthy-Li-Hodge-1} was to modify the definitions of k-forms by replacing the spaces $\wedge^k\mathcal{H}$ by Hilbert spaces $\mathcal{H}^{(k)}$, for $k=1,2,\dots$, continuously included in the projective exterior powers $\wedge^kTC_{x_0}$. For the `projective exterior powers' the completion is made using the largest cross norm and the usual, geometric, differential forms are sections of the dual bundles  $(\wedge^kTC_{x_0})^\ast$. Our $H$-k-forms will be sections of the dual bundles $\mathcal{H}^{(k)}{^*}$.

To define $\mathcal{H}^{(k)}$, for simplicity we will deal only with the case of Brownian motion measures and Levi-Civita connections. The more general situation is touched upon in \cite {Elworthy-Li-vector-fields}; for details of the following see \cite {Elworthy-Li-Hodge-2}. Take an SDE as in  section~\ref{sec-Ito.maps} with corresponding It\^o map $\mathcal{I}$. It is shown that the map $h_1\wedge\dots\wedge h_k \mapsto T_\omega\mathcal{I}(h_1)\wedge\dots\wedge T_\omega(h_k)$ determines a continuous linear map $\wedge^k T_\omega\mathcal{I}:\wedge^kH\rightarrow \wedge^kT_{\mathcal{I}(\omega)}C_{x_0}$ from Hilbert space to Banach space. As done in section \ref{sec-Ito.maps} integrate over the fibres of $\mathcal{I}$ to define $$\overline{\wedge T\mathcal{I}}_\sigma:\wedge^kH\rightarrow\wedge^kT_{\sigma} C_{x_0}$$ for almost all $\sigma\in C_{x_0}$, by the conditional expectation: $$\overline{\wedge T\mathcal{I}}_\sigma(\underline{h})=\mathbb{E}\{\wedge^kT\mathcal{I}(\underline{h})|\mathcal{I}=\sigma\}$$ for $\underline{h}\in \wedge^kH$.
We then let $\mathcal{H}^{(k)}_\sigma$ be the image of $ \overline{\wedge T\mathcal{I}}_\sigma$ with its quotient Hilbert space structure. Thus $\mathcal{H}^{(1)}=\mathcal{H}$. As with the case $k=1$ these spaces depend only on the Riemannian structure of $M$, not on the choice of SDE we used to construct them (provided $\breve{\nabla}$ is the Levi-Civita connection).

For $k=2$ there is a detailed description. For this let $\pathR:\wedge^2TC_{x_0}\to \mathbb{L}(\mathcal{H};\mathcal{H})$ be the curvature operator of the damped Markovian connection on $\mathcal{H}$, see Section~\ref{sec-Ito.maps}, and let $\mathbb{T}:\mathcal{H}\times\HH\to \HH$ be its torsion.
We have
$$ \mathcal{H}^{(2)}=\{U\in\wedge^2TC_{x_0}: U-\pathR(U)\in \wedge^2\mathcal{H}\}$$
with inner product having the norm $|U|_{\mathcal{H}^2}=|U-\pathR(U)|_{\wedge^2\mathcal{H}}$. 
Alternatively, inverting $\id-\pathR$, we have $$ \mathcal{H}^2=\{V+\mathbf{ Q}(V): V\in\wedge^2\mathcal{H}\}$$
where the linear map $\mathbf{ Q}$ can be expressed in terms of the curvature of $M$ and involves a `damped translation' of 2-vectors on $M$ where the damping is by the second Weitzenbock curvature, just as the first, the Ricci curvature, appears in equation (\ref{eq:damp}).
It turns out that $`\Div$ '$ \mathbf{ Q}(u\wedge v) =\frac{1}{2}\mathbb{T}(u,v)$ for any bounded \emph{ adapted } $H$-vector fields $u$ and $v$ in the sense that for any smooth cylindrical one-form $\phi$ on $C_{x_0}$ we have $$\int_{C_{x_0}} d\phi(\mathbf{ Q}(u\wedge v) d\mu_{x_0} =- \frac{1}{2}\int_{C_{x_0}}\phi(\mathbb{T}(u,v))d\mu_{x_0}.$$
This relates to a result of Cruzeiro-Fang, \cite {Cruzeiro-Fang-01}, that for suitable $u$ and $v$ the torsion $\mathbb{T}(u,v)$ has `divergence' zero, in the corresponding sense. 

If we define the exterior $H$-derivative as usual on cylindrical one forms $\phi$ but restrict the resulting $(d_H \phi)_\sigma : \wedge^2 T_\sigma C_{x_0}\rightarrow \R$ to $\mathcal{H}^{(2)}_\sigma$ we obtain a map, with domain the smooth cylindrical one forms, into the $L^2$ $H$ two-forms, $L^2\ {\mathcal{H}^{(2)}}^\ast$. The cylindrical one forms when restricted to $\mathcal{H}^{*}$ form a dense subspace of $L^2\mathcal{H}^\ast$ and it turns out that this map  is closable as an operator on $L^2\mathcal{H}^\ast$. We obtain a closed exterior derivative operator $$d^1: \Dom(d)\subset L^2\mathcal{H}^\ast \rightarrow L^2{\mathcal{H}^{(2)}}^\ast$$
with a dual operator $\Div:\Dom (\Div)\subset L^2\mathcal{H}^{(2)}\rightarrow L^2\mathcal{H}$.

The covariant derivative determined by the damped Markovian connection on $\mathcal{H}$ can be considered as a closed operator $\nabla$ from its domain, $\D ^{2,1}\mathcal{H}$,  in $ L^2\mathcal{H}$ to $L^2(\mathcal{H}\otimes\mathcal{H})$ and so has an adjoint $\nabla^\ast$. The following suggests that our construction is a natural one, but the condition of adaptedness on the vector fields is essential:
\begin {proposition} 
[\cite {Elworthy-Li-Hodge-2}]
Let $u$ and $v$ be bounded and adapted $H$-vector fields on $C_{x_0}$. Suppose $u,v \in \D ^{2,1}\mathcal{H}$ then $u\wedge v\in \Dom\nabla^\ast$ and 
$$\nabla^\ast(u\wedge v)=\Div \Big((Id+\mathbf{ Q})(u\wedge v)\Big).$$
\end{proposition}

It turns out that the exterior product $\phi^1\wedge\phi^2$ of two $H$-one-forms can be considered as an $H$-two-form in a consistent way.  Essentially this is because although an element in some $\mathcal{H}^{(2)}_\sigma$ is not in $\mathcal{H}_\sigma\otimes\mathcal{H}_\sigma$,
a space which can be identified with the Hilbert-Schmidt maps on $\mathcal{H}_\sigma$,
it can be identified with a bounded linear map on $\mathcal{H}_\sigma$, and elements of the uncompleted tensor product of $\mathcal{H}_\sigma$ with itself act as linear functionals on the bounded linear maps. We have then:
\begin{proposition} [\cite {Elworthy-Li-Hodge-2}]
Suppose $f\in\D ^{2,1}(C_{x_0};\R)$ and $\phi$ is a bounded $H$-one-form which is in the domain of the exterior derivative, and is bounded together with $d\phi$. Then $f\phi$ is in the domain of the exterior derivative and 
$$d^1(f\phi)=df\wedge \phi +fd^1\phi.$$
\end {proposition}

\subsection{A Hodge-Kodaira decomposition for one and two forms}

The key step to prove closability of the exterior derivative on these $H$-k-forms is to prove an analogue of Theorem \ref{th:TI}. We would like a rich set of $L^2$ maps $\underline{h}:C_{x_0}\to\wedge^k\HH$ such that $$\overline{\wedge^k T\I(\underline{h})}_\sigma:=\mathbb{E}\{\wedge^k T\I(\underline{h})|\I=\sigma\}\in\HH^{(k)}_\sigma$$
almost surely. For $k=1$ this holds for all such $\underline{h}$ by Theorem \ref{th:TI}. For $k=2$ it is claimed for an adequately rich family in \cite{Elworthy-Li-Hodge-2}, for all relevant It\^o maps, and for all $\underline{h}$ if the It\^o map is defined via a symmetric space structure. It is unknown for higher $k$ largely because of the apparently complicated algebraic structure of the spaces $\HH^{(k)}$ for higher $k$. (On the other hand in \cite{Elworthy-Li-vector-fields} it is shown that an important class of $k$-vector fields, defined for $k=1,..,n-1$ are $L^2$ sections of $\HH^{(k)}$, when $k=1,2$: these are important in the sense that they give integration by parts results, or generalised `Bismut-formulae',
for the finite dimensional exterior derivatives $dP_t\phi$ of the the heat semigroup on forms on $M$ in terms of a path integral of $\phi$ itself.)

From these results for $k=1,2$ we have now closed operators $$d^k:\Dom(d^k)\subset L^2\Gamma\HH^{(k)} \to L^2\Gamma\HH^{(k+1)}$$ for $k=1,2$ with $d^2d^1=0$. This leads to the Hodge-Kodaira decomposition: 
\begin{equation} 
\label{eq:Hodge-Kodaira}
L^2 \Gamma\HH^k=\overline{\image(d^{(k-1)})}\oplus \overline{\image((d^k)^\ast}\oplus (\ker d^k\cap \ker (d^{(k-1)})^\ast)
\end{equation}
for $ k=1,2$, as given for $k=1$ in \cite {Elworthy-Li-Hodge-1}, and for $k=2$ in \cite {Elworthy-Li-Hodge-2}. Here $d^0$ refers to $d$, and in the case $k=1$ the image of $d$ is closed by Fang's theorem, e.g. see the Clark-Ocone formula  in \cite
{Elworthy-LeJan-Li-book}. Moreover we have self-adjoint operators $(d^k)^\ast d^k+d^{(k-1)}(d^{(k-1)})^\ast$ acting on the spaces of $L^2$ $H$-k-forms for $k=1,2$. For $k=1$ the decomposition plus Fang's theorem shows that the space of  $L^2$ harmonic one-forms represents the  $L^2$ deRham cohomology group of $H$-one-forms.

\subsection{Lie groups with flat connection}
At present we have no information about even the first $L^2$ deRham group for non-flat manifolds. However in \cite {Fang-Franchi-paths}, Fang \& Franchi considered the case where $M$ is a compact Lie group $G$ with bi-invariant metric. For the Bismut tangent spaces coming from a right invariant flat connection the natural It\^o map to use is that of a left invariant SDE $dx_t=TL_{x_t}\circ dB_t $ for $B_.$ a Brownian motion on the Lie algebra   
$\g$. There is no `redundant noise' and the derivative of the It\^o map maps the Cameron-Martin space into the Bismut tangent spaces, and its exterior powers onto those of the Bismut spaces. There is no problem with the definition of the exterior derivative and they showed that the It\^o map can be used to transfer Shigekawa's results for classical Wiener space, determining a Hodge-Kodaira decomposition, and giving the vanishing of $L^2$ harmonic forms and so of the corresponding deRham cohomology groups.

\section{Loop spaces}
We have not extended the It\^o map techniques described above in any systematic way to the case of loop spaces, (but see \cite{Aida-irreducibility}), and here will only briefly describe the basic set up and some relevant results. The surveys \cite{Leandre-loops-survey} and \cite{Aida-loops-survey} give more information and references. Special motivation for the development of analysis on these spaces has come from the loop space approach to index theorems, as in \cite {Bismut-loops-85}, and the Hohn-Stolz conjecture, \cite{Stolz-conjecture}. However note that this theory is based on tangent spaces of vectors  which are in some sense in $L^{2,1}$ and it is not clear that this is always what is relevant to some physical or topological situations, e.g. see \cite{Freed}.

On the space of based loops, or more generally on the spaces $C_{{x_0},{y_0}}$ of continuous paths $\sigma:[0,T] \rightarrow M$ with $\sigma(0)=x_0$ and $\sigma(1)=y_0$, for $y_0\in M$, a natural measure to take is the \emph{Brownian Bridge measure}, $\mu_{{x_0},{y_0}}$, obtained by conditioning Brownian motion from $x_0$ to be at $y_0$ at time $T$. If equation (\ref{eq:SDE}) has solutions which are Brownian motions then the equation: \begin{equation} \label{eq:bridge} db_t=X(b_t)\circ dB_t+ A(b_t)dt+\bigtriangledown \log p_{T-t}(b_t;y_0)dt \end{equation} will have It\^o map which sends Wiener measure to $\mu_{{x_0},{y_0}}$. Here $p_t(x,y)$ is the heat kernel as in section \ref{sec-Bismut}. 

For the space $L(M)$ of free loops, i.e of continuous $\sigma: S^1\rightarrow M$, there is Bismut's measure, $\mu_L$, which can be defined as $\int_Mp_T(y,y)\mu_{{y},{y}}\,dy$ with $T=2\pi$, \cite{Bismut-loops-85}. This measure is invariant under the action of $S^1$.  A variant of this when $M$ is a Lie group is to average using normalised Haar measure rather than the heat kernel. 
Either of these loop spaces could be furnished with a \emph{heat kernel measure}, $\mu_h$. This is defined by choosing a base point, e.g. a constant loop, and constructing a `Brownian motion' on the loop space starting at that point, running it for some fixed time, $\tau$ say, and using its probability distribution as $\mu_h$. This will depend on $\tau$ and for free loops an extra averaging over the initial base point to retain $S^1$-invariance is needed. The construction of such Brownian motions goes back to Baxendale, see \cite {Baxendale-Edinburgh} and Gaveau \& Mazet, \cite{Gaveau-Mazet} but has been most developed for loop groups, \cite{Malliavin-loops}. For based paths on a compact simply connected Lie group this measure has been shown to be equivalent to the Brownian Bridge measure, \cite {Aida-Driver}.

In these contexts, by results of several people including M. P. \& P. Malliavin, Driver, Hsu, Leandre, Enchev \& Stroock, and Aida  there are integration by parts formulae and associated Sobolev spaces based on Bismut tangent spaces defined similarly to those above, though for  Lie groups flat connections are often used to define the Bismut tangent spaces. For the Brownian bridge and Bismut measures there are cohomology results: using stochastic Chen forms in \cite{Jones-Leandre}, `Sobolev differential forms'  \cite {Leandre-96}, and more recently `Chen-Souriau cohomology' defined via a `stochastic diffeology', see \cite{Leandre-loops-survey}. In general the resulting cohomology agrees with the usual singular real cohomology. We refer the reader to MathScinet to see the variety of constructions in these and related situations by R. Leandre.

 For compact Lie groups with bi-invariant metrics, and with Bismut tangent spaces defined by flat left, or right, invariant connections, Fang\& Franchi were able to extend their results for path spaces to based loops defining  Hodge-Kodaira operators on forms 
and giving a `Weitzenbock formulae' for them, \cite {Fang-Franchi-loops}. This formula rather clearly shows the form of these operators as `Witten' or `Bismut'  Laplacians where the `perturbing' vector field is not an $H$-vector field, and gives rise to stochastic integrals in the formulae. The curvature part of the formulae has a Ricci term, which requires a careful summation, as in \cite{Freed}.

One striking result for the Brownian bridge measure is the following by Eberle:
\begin{theorem}[\cite{Eberle-gap}] \label{th:eberle} 
Suppose the compact manifold $M$ has a closed geodesic for which there is a neighbourhood in $M$ of constant negative curvature. Then on the loop spaces $C_{{x_0},{x_0}}$ with Brownian bridge measure, and $L(M)$ with Bismut measure, the self-adjoint operator $\CL=-d^*d$ does not have a spectral gap.
\end {theorem}
Spectral gaps for the Hodge-Kodaira `Laplacians' are important in Hodge theory since they correspond to the (exterior) derivative operators having closed range in $L^2$. At present it is unknown if there is ever a spectral gap for $\CL$ for these measures for loops on non-flat manifolds, e.g. on spheres. However for heat kernel measures on compact Lie groups with bi-invariant metrics Driver \& Lohrenz proved the existence of a Log Sobolev inequality and so of a spectral gap for $\CL$, see \cite{Fang-loops}.

An alternative approach to based loops has been to represent them by `submanifolds' of
classical Wiener space by  choosing a suitable (i.e. a quasi-continuous) version of the stochastic development and considering the inverse image $\widetilde{C}$, say,  under it, of the based loops on $M$. This construction is dependent on `quasi-sure' analysis, see \cite{Malliavin-book}, where our measure theoretic concepts are refined potential theoretically, so that $\widetilde{C}$ can be defined up to sets of capacity zero. To a certain extent this allows $\widetilde{C}$ to be treated as a submanifold of co-dimension the dimension of $M$, with a differential form theory and Weitzenbock formula, see \cite {Kazumi-Shigekawa}, and \cite{Kusuoka-ICM}.

\end {document}